\documentclass{monsky2009} 

\usepackage{graphicx}

\usepackage{amssymb}
\usepackage{amsmath}

\usepackage{bm}

\newenvironment{theorem*}[1]{\textbf{#1}\itshape \hspace{.3em}}{\upshape}
\newenvironment{remark*}[1]{\textbf{#1}\itshape \hspace{.3em}}{\upshape}
\newenvironment{corollary*}[1]{\textbf{#1}\itshape \hspace{.3em}}{\upshape}
\newenvironment{proof}{\textbf{Proof\hspace{.3em}}}{}

\newtheorem{definition}{Definition}[section]
\newtheorem{theorem}[definition]{Theorem}
\newtheorem{lemma}[definition]{Lemma}

\newtheorem{example}[definition]{Example}

\newcommand{\conste}{\ensuremath{\mathrm{e}}}
\newcommand{\consti}{\ensuremath{\mathrm{i}}}

\newcommand{\newo}{\ensuremath{\mathcal{O}}}

\newcommand{\newc}{\ensuremath{\mathbb{C}}}
\newcommand{\newz}{\ensuremath{\mathbb{Z}}}

\begin{document}

\begin{frontmatter}






\title{The limit as $p\rightarrow\infty$ of the Hilbert-Kunz multiplicity of $\sum x_{i}^{d_{i}}$}
\author{Ira M. Gessel and Paul Monsky}

\address{Brandeis University, Waltham MA  02454-9110, USA\\ gessel@brandeis.edu and monsky@brandeis.edu }

\begin{abstract}
Let $p$ be a prime. The Hilbert-Kunz multiplicity, $\mu$, of the element $\sum x_{i}^{d_{i}}$ of \linebreak $({\newz}/{p})[x_{1},\ldots, x_{s}]$ depends on $p$ in a complicated way. We calculate the limit of $\mu$ as $p\rightarrow\infty$. In particular when each $d_{i}$ is 2 we show that the limit is 1 + the coefficient of $z^{s-1}$ in the power series expansion of $\sec z + \tan z$.
\end{abstract}


\end{frontmatter}


\section{Introduction}
\label{introduction}

Suppose $s\ge 2$, $d_{1},\ldots , d_{s}$ are positive integers, and $h$ is the element $\sum x_{i}^{d_{i}}$ of $A=({\newz}/{p})[x_{1},\ldots, x_{s}]$. Let $e_{n}(h)$ be the colength of the ideal generated by $h$ and the $x_{i}^{q}$ where $q=p^{n}$. For fixed $p$, Hilbert-Kunz theory tells us that $e_{n}=\mu q^{s-1}+O(q^{s-2})$ for some $\mu >0$; $\mu$ is the Hilbert-Kunz multiplicity of $h$. When $s=2$, $\mu = \min (d_{1},d_{2})$ and so is independent of $p$, but the dependence on $p$ is subtle when $s\ge 3$.

In her thesis, Han calculated $\mu$ (and in fact all of the $e_{n}$) when $s=3$. This result was extended to $s>3$ in \cite{3}. The second author realized afterwards that a result from \cite{3} gives an easy proof that $\mu\rightarrow$ a limit as $p\rightarrow\infty$, and a formula for the limit. The formula has been discovered by others since, but perhaps because their arguments were more complicated, they haven't presented them for publication.

Much of the interest of the above result lies in an elegant expression for the limit when each $d_{i}$ is 2; the limit is 1 + the coefficient of $z$ in the power series expansion of $\sec z + \tan z$. This was conjectured by the second author, and the first used Eulerian polynomials to provide a proof. At the request of several colleagues we're here writing down our old proofs. The limits of Hilbert-Kunz multiplicities have been studied in other situations; see Trivedi \cite{4}.

We begin with some easy results. When any $d_{i}$ is 1, $e_{n}=q^{s-1}$ irrespective of $p$, so the limit of $\mu$ is 1. Assume from now on that each $d_{i}>1$. The following function was studied in \cite{3}.

\begin{definition}
\label{def1.1}
$D_{p}(k_{1},\ldots, k_{s})=\mathrm{length}\ A/(\sum x_{i},x_{1}^{k_{i}},\ldots, x_{s}^{k_{s}})$.
\end{definition}

Evidently, $D_{p}(k_{1},\ldots, k_{s})$ is the length of $$({\newz}/{p})[x_{2},\ldots, x_{s}]{\displaystyle /}\left((x_{2}+\cdots + x_{s})^{k_{1}},x_{2}^{k_{2}},\ldots, x_{s}^{k_{s}}\right).$$ Since $({\newz}/{p})[x_{2},\ldots, x_{s}]$ is a free module of rank $p^{s-1}$ over $({\newz}/{p})[x_{2}^{p},\ldots, x_{s}^{p}]$, it follows that $D_{p}(pk_{1},\ldots, pk_{s}) = p^{s-1}D_{p}(k_{1},\ldots, k_{s})$.

\begin{lemma}
\label{lemma1.2}
Suppose $d_{i}u_{i}\le p\le d_{i}v_{i}$. Then for $n>0$, each of the $e_{n}(h)/q^{s-1}$ lies between $dp^{1-s}D_{p}(u_{1},\ldots , u_{s})$ and $dp^{1-s}D_{p}(v_{1},\ldots , v_{s})$, where $d$ is the product of the $d_{i}$.
\end{lemma}

\begin{proof}
For $n>0$, $e_{n}$ is bounded below by the colength of the ideal generated by $h$ and the $x_{i}^{\frac{qd_{i}u_{i}}{p}}$. Since $A$ is free of rank $d$ over $({\newz}/{p})[x_{1}^{d_{1}},\ldots ,x_{s}^{d_{s}}]$ this colength is $dD_{p}\left(\frac{qu_{1}}{p}, \cdots, \frac{qu_{s}}{p}\right) = \left(\frac{q}{p}\right)^{s-1}dD_{p}(u_{1},\ldots, u_{s})$. Dividing by $q^{s-1}$, we get the lower bound, and the upper bound is derived similarly.
\qed
\end{proof}

Now the $A$-module $\left(\sum x_{i}, x_{1}^{k_{1}},\ldots, x_{s}^{k_{s}}\right)/\left(\sum x_{i},x_{1}^{k_{1}+1},x_{2}^{k_{2}},\ldots, x_{s}^{k_{s}}\right)$ is annihilated by $x_{1}$, and so may be viewed as a $({\newz}/{p})[x_{2},\ldots, x_{s}]$-module. As such, it is cyclic, generated by $x_{1}^{k_{1}}$ and annihilated by $x_{2}+\cdots + x_{s}$ and by each of $x_{2}^{k_{2}},\ldots, x_{s}^{k_{s}}$. So its length is at most $D_{p}(k_{2},\ldots, k_{s}).$

\begin{lemma}
\label{lemma1.3}
If each of $u_{1},\ldots, u_{s}$ is $<p$ then $D_{p}(u_{1}+1, u_{2}+1, \ldots, u_{s}+1)-D_{p}(u_{1}, \ldots, u_{s})\le sp^{s-2}$.

\end{lemma}

\begin{proof}
The argument preceding the lemma shows that $D_{p}(u_{1}+1,u_{2},\ldots , u_{s})-D_{p}(u_{1},u_{2},\ldots , u_{s})\le D_{p}(u_{2},\ldots, u_{s})$. Since each $u_{i}$ is $<p$, this is $\le p^{s-2}$. Combining this with $s-1$ similar inequalities we get the result.
\qed
\end{proof}

\begin{theorem}
\label{theorem1.4}
$\frac{e_{1}(h)}{p^{s-1}}$ and $\mu$ differ by at most $\frac{ds}{p}$. Consequently, $\lim_{p\rightarrow \infty}(\mu)=\lim_{p\rightarrow\infty}\left(\frac{e_{1}(h)}{p^{s-1}}\right)$, provided the latter limit exists.

\end{theorem}

\begin{proof}
Set $u_{i}=\lfloor\frac{p}{d_{i}}\rfloor$. Since $d_{i}>1$, each $u_{i}<p$. If we let $v_{i}$ be $u_{i}+1$, then Lemmas \ref{lemma1.2} and \ref{lemma1.3} show that all the $e_{n}(h)/q^{s-1}$, $n>0$, lie in an interval of length $\le (sp^{s-2})(dp^{1-s})=\frac{ds}{p}$. But $\mu$ is in the closure of this interval.
\qed
\end{proof}

\begin{theorem}
\label{theorem1.5}
Suppose for each $p$ we are given integers $a_{1},\ldots, a_{s}\ge 0$ with $a_{i}=\frac{p}{d_{i}}+O(1)$. Then there is a $k$ such that for all $p$ the difference between $\mu$ and $dp^{1-s}D_{p}(a_{1},\ldots, a_{s})$ is at most $\frac{k}{p}$.
\end{theorem}

\begin{proof}
Fix $N$ large and let $u_{i}=a_{i}-N$, $v_{i}=a_{i}+N$. Then $d_{i}u_{i}\le p\le d_{i}v_{i}$. And when $p$ is large, $u_{i}\ge 0$ and $v_{i}\le p$. Now $dp^{1-s}D_{p}(a_{1},\ldots, a_{s})$ and each $e_{n}/q^{s-1}$, $n>0$, lie between $dp^{1-s}D_{p}(u_{1},\ldots, u_{s})$ and $dp^{1-s}D_{p}(v_{1},\ldots, v_{s})$. The argument of Theorem \ref{theorem1.4} shows they lie in an interval of length $\le 2N\left(\frac{ds}{p}\right)$. The closure of this interval contains $\mu$.
\qed
\end{proof}

When $a_{1},\ldots, a_{s}$ are $\le p$, Theorem 2.20 of \cite{3} gives the following formula for $D_{p}(a_{1},\ldots,a_{s})$. Let $\gamma$ be $\lfloor\frac{1}{2}\sum(a_{i}-1)\rfloor$. Then $D_{p}(a_{1},\ldots, a_{s})$ is the sum as $\lambda$ runs over $\newz$ of the coefficients of the $t^{\gamma +\lambda p}$ in the polynomial $\Pi\left(\frac{1-t^{a_{i}}}{1-t}\right)$. In the next section we'll combine this result with Theorem \ref{theorem1.5} to calculate the limit of $\mu$ as $p\rightarrow\infty$.

\section{The limit formula}
\label{section2}

\begin{definition}
\label{def2.1}
For $\lambda$ in $\newz$, $C_{\lambda}=\sum (\epsilon_{1}\cdots \epsilon_{s})\left(\frac{\epsilon_{1}}{d_{1}}+\cdots + \frac{\epsilon_{s}}{d_{s}}-2\lambda\right)^{s-1}$, where the sum extends over the $s$-tuples $\epsilon_{1},\ldots, \epsilon_{s}$ with each $\epsilon_{i}$ in $\{-1,1\}$ and $\frac{\epsilon_{1}}{d_{1}}+\cdots + \frac{\epsilon_{s}}{d_{s}}>2\lambda$.

\end{definition}

For each $p$ choose integers $a_{1},\ldots, a_{s}\ge 0$ so that $\sum a_{i}\equiv s\pod{2}$, and $a_{i}=\frac{p}{d_{i}}+O(1)$. Suppose $\epsilon_{1},\ldots, \epsilon_{s}$ are in $\{-1,1\}$ and $\lambda$ is in $\newz$. Let $a$ in $Q$ be $\frac{\epsilon_{1}}{d_{1}}+\cdots + \frac{\epsilon_{s}}{d_{s}}-2\lambda$; this is independent of $p$. Let $\alpha$ be $\left(\frac{1}{2}\sum(\epsilon_{i}a_{i}-1)\right)-p\lambda$. Since $\sum a_{i}\equiv s\pod{2}$, $\alpha$ is in $\newz$. Evidently $\alpha =\frac{pa}{2}+O(1)$.

We fix $\epsilon_{1}\cdots \epsilon_{s}$ and $\lambda$, and study how the coefficient of $t^{\alpha}$ in $(1-t)^{-s}$ depends on $p$. When $a<0$, $\alpha<0$ for large $p$ and the coefficient is 0. When $a=0$, $\alpha$ is $O(1)$ and the coefficient is $O(1)$; since $s\ge 2$ it is $O(p^{s-2})$. Now suppose $a>0$. Then for large $p$, $\alpha > 0$ and the coefficient is $\binom{\alpha + s-1}{s-1}$. Since $\alpha +s-1$ is $\frac{pa}{2}+O(1)$, we get $\frac{1}{(s-1)!}\left(\frac{pa}{2}\right)^{s-1}+O(p^{s-2})=\frac{2^{1-s}}{(s-1)!}a^{s-1}p^{s-1}+O(p^{s-2})$.

\begin{lemma}
\label{lemma2.2}
Let $\gamma = \frac{1}{2}\sum(a_{i}-1)$. Then the coefficient of $t^{\gamma}$ in $\Pi\left(\frac{1-t^{a_{i}}}{1-t}\right)$ is $\frac{2^{1-s}}{(s-1)!}C_{0}p^{s-1}+O(p^{s-2})$ with $C_{0}$ as in Definition \ref{def2.1}.

\end{lemma}

\begin{proof}
$\Pi\left(\frac{1-t^{a_{i}}}{1-t}\right)=(1-t)^{-s}\Pi(1-t^{a_{i}})$. Multiplying the second product out we express our coefficient in terms of coefficients of $(1-t)^{-s}$. Explicitly it is (the coefficient of $t^{\gamma}$ in $(1-t)^{-s}$)$-$(the sum of the coefficients of the $t^{\gamma-a_{i}}$)$+$(the sum of the coefficients of the $t^{\gamma-a_{i}-a_{j}}$)$-\cdots$. The paragraph before the lemma, with $\lambda=0$, tells us the behavior of each term as $p\rightarrow\infty$ and gives the result.
\qed
\end{proof}

More generally:

\begin{lemma}
\label{lemma2.3}
The coefficient of $t^{\gamma-\lambda p}$ in $\Pi\left(\frac{1-t^{a_{i}}}{1-t}\right)$ is $\frac{2^{1-s}}{(s-1)!}C_{\lambda}p^{s-1}+O(p^{s-2})$. Furthermore $C_{\lambda}=C_{-\lambda}$.
\end{lemma}

\begin{proof}
The argument of Lemma \ref{lemma2.2} gives the first result. Since the coefficients of $t^{\gamma+N}$ and $t^{\gamma-N}$ in $\Pi\left(\frac{1-t^{a_{i}}}{1-t}\right)$ are equal, the second result follows.
\qed
\end{proof}

Note that when $\lambda \ge \frac{s}{4}$, $\frac{\epsilon_{1}}{d_{1}}+\cdots + \frac{\epsilon_{s}}{d_{s}}\le \frac{s}{2}\le 2\lambda$ and so $C_{\lambda}=0$. By Lemma \ref{lemma2.3}, $C_{\lambda}=0$ when $|\lambda|\ge\frac{s}{4}$.

\begin{theorem}
\label{theorem2.4}
As $p\rightarrow\infty$, $\mu\rightarrow\frac{d(2^{1-s})}{(s-1)!}\left(\sum C_{\lambda}\right)=\frac{d(2^{1-s})}{(s-1)!}\left(C_{0}+2\sum_{\lambda>0}C_{\lambda}\right)$.
\end{theorem}

\begin{proof}
Since $C_{\lambda}=C_{-\lambda}$, the sums are equal. By Theorem \ref{theorem1.5} it suffices to show that $p^{1-s}D_{p}(a_{1},\ldots, a_{s})\rightarrow \frac{2^{1-s}}{(s-1)!}\sum C_{\lambda}$ as $p\rightarrow\infty$. But this follows immediately from Lemmas \ref{lemma2.2}, \ref{lemma2.3} and the result from \cite{3} quoted at the end of the Introduction.
\qed
\end{proof}

\begin{example}
\label{example2.5}
Supppose $s=4$ and each $d_{i}$ is $4$. Then $C_{0}=\left(\frac{1}{4}+\frac{1}{4}+\frac{1}{4}+\frac{1}{4}\right)^{3}-4\left(\frac{1}{4}+\frac{1}{4}+\frac{1}{4}-\frac{1}{4}\right)^{3}=\frac{1}{2}$, while $C_{\lambda}=0$ for $\lambda\ne 0$. So $\mu\rightarrow\frac{256}{8\cdot 6}\cdot\frac{1}{2}=\frac{8}{3}$. In fact, $\mu =\frac{8}{3}\left(\frac{2p^{2}+2p+3}{2p^{2}+2p+1}\right)$ if $p\equiv 1\pod{4}$ and $\frac{8}{3}\left(\frac{2p^{2}-2p+3}{2p^{2}-2p+1}\right)$ if $p\equiv 3\pod{4}$.
\end{example}

From now on we assume each $d_{i}$ is $2$.

\begin{definition}
\label{def2.6}
If $a$ is an integer, $f_{s}(a)=a^{s-1}-\binom{s}{1}(a-2)^{s-1}+\binom{s}{2}(a-4)^{s-1}-\cdots$, where we make the convention that $c^{s-1}=0$ when $c<0$.
\end{definition}

\begin{theorem}
\label{theorem2.7}
As $p\rightarrow\infty$, $$\mu\rightarrow\frac{1}{(s-1)!}\cdot\frac{1}{2^{s-2}}\left(f_{s}(s)+2f_{s}(s-4)+2f_{s}(s-8)+2f_{s}(s-12)+\cdots\right).$$ This may also be written as $\frac{1}{(s-1)!}\cdot\frac{1}{2^{s-2}}\cdot\sum f_{s}(a)$, the sum extending over all $a\equiv s\pod{4}$.
\end{theorem}

\begin{proof}
$C_{0}=\left(\frac{s}{2}\right)^{s-1}-\binom{s}{1}\left(\frac{s-2}{2}\right)^{s-1}+\binom{s}{2}\left(\frac{s-4}{2}\right)^{s-1}-\cdots = \frac{1}{2^{s-1}}\cdot f_{s}(s)$; similarly $C_{\lambda}= \frac{1}{2^{s-1}}\cdot f_{s}(s-4\lambda)$. Now apply Theorem \ref{theorem2.4}, noting that $\frac{d\cdot 2^{1-s}}{(s-1)!}=\frac{2}{(s-1)!}$.
\qed
\end{proof}

\begin{example}
\label{example2.8}
Supppose $s=5$. $f_{5}(5)=5^{4}-5\cdot 3^{4}+10\cdot 1^{4}=230$, while $f_{5}(1)=1^{4}=1$. So by Theorem \ref{theorem2.7}, $\mu\rightarrow\frac{1}{24}\cdot\frac{1}{8}(232)=\frac{29}{24}$. In fact, if $p>2$, $\mu=\frac{29p^{2}+15}{24p^{2}+12}$.
\end{example}

Proceeding as in Example \ref{example2.8}, the second author calculated the limit of $\mu$ for each $s\le 10$, finding that in each case the limit was $1+$ the coefficient of $z^{s-1}$ in the power series expansion of $\sec z + \tan z$. In the next section we'll use Eulerian polynomials to show that this holds for all $s$; this insight is due to the first author.

\section{The case $h=\sum x_{i}^{2}$}

\begin{definition}
\label{def3.1}
For $n\ge 0$, $A_{n}=(1-T)^{n+1}(1 + 2^{n}T+3^{n}T^{2}+\cdots)$.
\end{definition}

For example, $A_{4}=1+11T+11T^{2}+T^{3}$.

\begin{lemma}
\label{lemma3.2}
$A_{n}$ is a polynomial, and $A_{n}(1)=n!$.
\end{lemma}

\begin{proof}
Let $\Delta$ be the operator $f\rightarrow f(T+1)-f(T)$ on $\newz[T]$. The $n$-fold iterate of $\Delta$ evidently takes $T^{n}$ to the constant $n!$. It follows that $(1-T)^{n}(1+2^{n}T+3^{n}T^{2}+\cdots)=(\mbox{a polynomial in $T$})+\frac{n!}{(1-T)}$. Multiplying by $1-T$ and evaluating at $T=1$ we get the result.
\qed
\end{proof}

Euler \cite{1}, \cite{2} evaluated these Eulerian polynomials at $-1$. The values at $\consti$ are less familiar but we'll show how to derive them by an easy method. As divergent series are out of fashion, we'll proceed formally. Let $\newo$ be the complete local ring $\newc[[T,z]]$. If $u$ is in the maximal ideal of $\newo$, $\conste^{u}$ will denote the element $\sum_{n\ge 0}\frac{u^{n}}{n!}$ of $\newo$.

\begin{lemma}
\label{lemma3.3}
In $\newo$, $\left(\sum_{n\ge 1}\frac{A_{n}(T)}{(1-T)^{n}}\frac{z^{n}}{n!}\right)(1-T\conste^{z})=\conste^{z}-1$.
\end{lemma}

\begin{proof}
$\sum_{n\ge 0}\frac{A_{n}(T)}{(1-T)^{n+1}}\frac{z^{n}}{n!}=\sum_{n\ge 0}\frac{z^{n}}{n!}\left(1^{n}+2^{n}T+3^{n}T^{2}+\cdots\right)=\conste^{z}+T\conste^{2z}+T^{2}\conste^{3z}+\cdots$.  Multiplying by $(1-T\conste^{z})(1-T)$ we find that $\left(\sum_{n\ge 0}\frac{A_{n}(T)}{(1-T)^{n}}\frac{z^{n}}{n!}\right)(1-T\conste^{z})=\conste^{z}(1-T)$. Subtracting off the $n=0$ term we get $\conste^{z}(1-T)-(1-T\conste^{z})=\conste^{z}-1$.
\qed
\end{proof}

\begin{lemma}
\label{lemma3.4}
In $\newc[[z]]$, $\sum_{n\ge 1}\frac{A_{n}(\consti)}{(1+\consti)^{n}}\cdot\frac{z^{n}}{n!}=\frac{1-\conste^{\consti z}}{\conste^{\consti z}-\consti}$.
\end{lemma}

\begin{proof}
There is a continuous ring automorphism of $\newo$ taking $T$ to $T$ and $z$ to $z(1-T)$. Applying this to Lemma \ref{lemma3.3} we find that $\left(\sum_{n\ge 1}A_{n}(T)\frac{z^{n}}{n!}\right)(1-T\conste^{z(1-T)})=\conste^{z(1-T)}-1$. Now this is an identity in $\newc[T][[z]]$. Applying the continuous ring homomorphism $\newc[T][[z]]\rightarrow\newc[[z]]$ that takes $T$ to $\consti$ and $z$ to $\frac{z}{1+\consti}$ we get the result.
\qed
\end{proof}

\begin{theorem}
\label{theorem3.5}
If $s\ge 2$, $\frac{A_{s-1}(\consti)}{(1+\consti)^{s-2}}\cdot\frac{1}{(s-1)!}$ is the coefficient of $z^{s-1}$ in the power series expansion of $\sec z+\tan z$.
\end{theorem}

\begin{proof}
Multiplying both sides of Lemma \ref{lemma3.4} by $1+\consti$ and adding $1$ we find that in $\newc[[z]]$, $1+\sum_{n\ge 1}\frac{A_{n}(\consti)}{(1+\consti)^{n-1}}\frac{z^{n}}{n!}=\frac{1-\consti\conste^{\consti z}}{\conste^{\consti z}-\consti}$. If by $\sin z$, $\cos z$, $\sec z$, $\tan z$ we mean the Taylor series expansions of these functions, then $\frac{1-\consti\conste^{\consti z}}{\conste^{\consti z}-1}=\frac{(1+\sin z)-\consti\cos z}{\cos z-\consti(1-\sin z)}$.

Since $\frac{1+\sin z}{\cos z}$ and $\frac{\cos z}{1-\sin z}$ are each $\sec z + \tan z$, $1+\sum_{1}^{\infty}\frac{A_{n}(\consti)}{(1+\consti)^{n-1}}\frac{z^{n}}{n!}=\sec z + \tan z$ in $\newc[[z]]$, and we compare the coefficients of $z^{s-1}$.
\qed
\end{proof}

\begin{lemma}
\label{lemma3.6}
$\sum_{a}f_{s}(a)T^{a-1}=(1+T)^{s}A_{s-1}$, where $f_{s}(a)$ is as in Definition \ref{def2.6}.
\end{lemma}

\begin{proof}
$f_{s}(a)=a^{s-1}-\binom{s}{1}(a-2)^{s-1}+\binom{s}{2}(a-4)^{s-1}-\cdots$. So $\sum f_{s}(a)T^{a-1}=(1^{s-1}+2^{s-1}T+3^{s-1}T^{2}+\cdots)-\binom{s}{1}(1^{s-1}T^{2}+2^{s-1}T^{3}+3^{s-1}T^{4}+\cdots)+\binom{s}{2}(1^{s-1}T^{4}+2^{s-1}T^{5}+3^{s-1}T^{6}+\cdots)-\cdots = (1^{s-1}+2^{s-1}T+3^{s-1}T^{2}+\cdots)(1-T^{2})^{s}=(1+T)^{s}A_{s-1}$.
\qed
\end{proof}

\begin{theorem}
\label{theorem3.7}
Let $c=\frac{A_{s-1}(\consti)}{(1+\consti)^{s-2}}$. Then $\sum f_{s}(a)$, the sum extending over all $a\equiv s\pod{4}$, is $2^{s-2}\left(A_{s-1}(1)+\frac{c}{2}+\frac{\bar{c}}{2}\right)$.
\end{theorem}

\begin{proof}
Let $P=\sum_{j}f_{s}(j-3s)T^{j}$. Since $4$ divides $j$ if and only if $j-3s\equiv s\pod{4}$, our sum is $\frac{1}{4}\left(P(1)+P(-1)+P(\consti)+P(-\consti)\right)$. Now $P=\sum f_{s}(j)T^{j+3s}$ which is $T^{3s+1}(1+T)^{s}A_{s-1}$ by Lemma \ref{lemma3.6}. Thus $\frac{1}{4}P(1)=2^{s-2}A_{s-1}(1)$ and $\frac{1}{4}P(-1)=0$. Furthermore, $P(\consti)=(-\consti)^{s-1}(1+\consti)^{s}A_{s-1}(\consti)$. Since $(-\consti)^{s-1}(1+\consti)^{s}(1+\consti)^{s-2}=(-\consti)^{s-1}(2\consti)^{s-1}=2^{s-1}$, $\frac{1}{4}P(\consti)=\frac{1}{4}\cdot 2^{s-1}\cdot\frac{A_{s-1}(\consti)}{(1+\consti)^{s-2}}=2^{s-2}\cdot\left(\frac{c}{2}\right)$. Conjugating we find that $\frac{1}{4}P(-\consti)=2^{s-2}\cdot\left(\frac{\bar{c}}{2}\right)$.
\qed
\end{proof}

\begin{theorem}
\label{theorem3.8}
Suppose $h=\sum_{1}^{s}x_{i}^{2}$. Then as $p\rightarrow\infty$, the Hilbert-Kunz multiplicity of $h\rightarrow 1 +$ the coefficient of $z^{s-1}$ in $\sec z + \tan z$.
\end{theorem}

\begin{proof}
Theorems \ref{theorem2.7} and \ref{theorem3.7} show that $\mu\rightarrow\frac{1}{(s-1)!}\left(A_{s-1}(1)+\frac{c}{2}+\frac{\bar{c}}{2}\right)$. But $A_{s-1}(1)=(s-1)!$. And Theorem \ref{theorem3.5} shows that $\frac{c}{(s-1)!}$ and $\frac{\bar{c}}{(s-1)!}$ are each the coefficient of $z^{s-1}$ in $\sec z + \tan z$.
\qed
\end{proof}



\label{}



\end{document}